\documentclass[11pt,reqno]{amsart}
\usepackage[utf8]{inputenc}
\usepackage[T1]{fontenc}
\usepackage{lmodern}
\usepackage{amsfonts,amsthm,amsmath,amssymb}
\usepackage{graphicx}
\usepackage{float}
\usepackage{fvextra}
\usepackage[dvipsnames,table]{xcolor}   
\usepackage{enumerate}
\usepackage{hyperref}
\usepackage{color}
\usepackage{tikz-cd}
\usetikzlibrary{calc}
\usepackage{transparent}
\usepackage[margin=1in]{geometry}
\usepackage[
    maxbibnames=99,
    backend=biber,
    style=alphabetic,
    sorting=nyt,
    giveninits=true
]{biblatex}
\DeclareFieldFormat{pages}{#1}
\renewbibmacro{in:}{%
  \ifentrytype{article}
    {}
    {\bibstring{in}%
     \printunit{\intitlepunct}}}
\DeclareFieldFormat
[article,inbook,incollection,inproceedings,patent,thesis,unpublished]
  {title}{\mkbibemph{#1}}
\DeclareFieldFormat{journaltitle}{#1\isdot}
\DeclareFieldFormat[article]{volume}{\mkbibbold{#1}}
\DeclareFieldFormat[article]{number}{\bibstring{number}\addnbspace #1}

\renewbibmacro*{journal+issuetitle}{%
  \usebibmacro{journal}%
  \setunit*{\addspace}%
  \iffieldundef{series}
    {}
    {\newunit
     \printfield{series}%
     \setunit{\addspace}}%
  \printfield{volume}%
  \setunit{\addspace}%
  \usebibmacro{issue+date}%
  \setunit{\addcomma\space}%
  \printfield{number}%
  \setunit{\addcolon\space}%
  \usebibmacro{issue}%
  \setunit{\addcomma\space}%
  \printfield{eid}
  \newunit}

\newtheoremstyle{mystyle}
  {}
  {}
  {\itshape}
  {}
  {\bfseries}
  {.}
  { }
  {\thmname{#1}\thmnumber{ #2}\thmnote{ (#3)}}

\theoremstyle{mystyle}
\newtheorem{Thm}{Theorem}[section]

\newtheorem{Conj}[Thm]{Conjecture}

\theoremstyle{definition}

\theoremstyle{remark}
\newtheorem{Rmk}[Thm]{Remark}

\newcommand{\Z}{\mathbb{Z}}
\newcommand{\wrap}{\operatorname{wrap}}

\tikzset{
  strand/.style={draw=black,line width=.85pt,line cap=round,line join=round},
  over/.style={strand,preaction={draw=white,line width=3.9pt,line cap=butt}},
  axis/.style={font=\Large,inner sep=0pt}
}

\newcommand{\Crossing}[4]{
  \pgfmathsetmacro{\xr}{#1+#3}
  \pgfmathsetmacro{\yt}{#2+.34}
  \pgfmathsetmacro{\yb}{#2-.34}
  \ifnum#4=1
    \draw[strand] (#1,\yt) .. controls (#1,#2+.12) and (\xr,#2-.12) .. (\xr,\yb);
    \draw[over] (#1,\yb) .. controls (#1,#2-.12) and (\xr,#2+.12) .. (\xr,\yt);
  \else
    \draw[strand] (#1,\yb) .. controls (#1,#2-.12) and (\xr,#2+.12) .. (\xr,\yt);
    \draw[over] (#1,\yt) .. controls (#1,#2+.12) and (\xr,#2-.12) .. (\xr,\yb);
  \fi
}

\newcommand{\SkeinCrossing}{%
  \begin{tikzpicture}[x=.45cm,y=.45cm,baseline=-.55ex]
    \path[use as bounding box] (-.52,-.52) rectangle (.52,.52);
    \draw[strand] (-.48,-.48)--(.48,.48);
    \draw[over] (-.48,.48)--(.48,-.48);
  \end{tikzpicture}%
}

\newcommand{\SkeinZero}{%
  \begin{tikzpicture}[x=.45cm,y=.45cm,baseline=-.55ex]
    \path[use as bounding box] (-.52,-.52) rectangle (.52,.52);
    \draw[strand] (-.48,-.48)..controls(-.12,-.24)and(-.12,.24)..(-.48,.48);
    \draw[strand] (.48,-.48)..controls(.12,-.24)and(.12,.24)..(.48,.48);
  \end{tikzpicture}%
}

\newcommand{\SkeinOne}{%
  \begin{tikzpicture}[x=.45cm,y=.45cm,baseline=-.55ex]
    \path[use as bounding box] (-.52,-.52) rectangle (.52,.52);
    \draw[strand] (-.48,-.48)..controls(-.24,-.12)and(.24,-.12)..(.48,-.48);
    \draw[strand] (-.48,.48)..controls(-.24,.12)and(.24,.12)..(.48,.48);
  \end{tikzpicture}%
}

\newcommand{\SkeinCircle}{%
  \begin{tikzpicture}[x=.45cm,y=.45cm,baseline=-.55ex]
    \path[use as bounding box] (-.17,-.48) rectangle (.72,.48);
    \draw[strand] (.275,0) circle[radius=.40];
  \end{tikzpicture}%
}

\newcommand{\FiveBlock}[2]{
  \begin{scope}[shift={(#1,0)}]
  \foreach \level/\y in {1/1.6,2/.8,3/0,4/-.8,5/-1.6}{
    \ifnum#2=1
      \ifnum\level=1 \def\pair{0}\def\sgn{1}\else
      \ifnum\level=5 \def\pair{0}\def\sgn{1}\else
        \def\pair{.75}\def\sgn{-1}\fi\fi
    \else
      \ifnum\level=1 \def\pair{.75}\def\sgn{-1}\else
      \ifnum\level=5 \def\pair{.75}\def\sgn{-1}\else
        \def\pair{0}\def\sgn{1}\fi\fi
    \fi
    \foreach \x in {0,.75,1.5}{
      \draw[strand] (\x,\y+.4)--(\x,\y+.34);
      \draw[strand] (\x,\y-.34)--(\x,\y-.4);
      \pgfmathtruncatemacro{\passive}{abs(\x-\pair)>.01 && abs(\x-\pair-.75)>.01}
      \ifnum\passive=1 \draw[strand] (\x,\y-.34)--(\x,\y+.34);\fi
    }
    \Crossing{\pair}{\y}{.75}{\sgn}
  }
  \end{scope}
}

\newcommand{\AnnularArc}[3]{
  \pgfmathsetmacro{\midpoint}{(#1+#3)/2}
  \pgfmathsetmacro{\handle}{.28*(#3-#1)}
  \draw[strand] (#1,3.2)
    .. controls (#1,#2) and (\midpoint-\handle,#2) .. (\midpoint,#2)
    .. controls (\midpoint+\handle,#2) and (#3,.55*#2) .. (#3,0)
    .. controls (#3,-.55*#2) and (\midpoint+\handle,-#2) .. (\midpoint,-#2)
    .. controls (\midpoint-\handle,-#2) and (#1,-#2) .. (#1,-3.2);
}

\newcommand{\AnnularLink}[1]{
  \begin{tikzpicture}[x=.77cm,y=.52cm]
    \FiveBlock{0}{0}
    \FiveBlock{2.25}{1}
    \draw[strand] (1.5,2)
      .. controls (1.5,2.65) and (2.25,2.65) .. (2.25,2);
    \draw[strand] (1.5,-2)
      .. controls (1.5,-2.65) and (2.25,-2.65) .. (2.25,-2);
    \foreach \x in {0,.75,3,3.75}{
      \draw[strand] (\x,-2)--(\x,-3.2);
    }
    \ifnum#1=0
      \foreach \x in {0,.75,3,3.75}{\draw[strand] (\x,2)--(\x,3.2);}
    \else
      \foreach \x in {0,.75,3,3.75}{
        \draw[strand] (\x,2)--(\x,2.46);
        \draw[strand] (\x,3.14)--(\x,3.2);
      }
      \Crossing{0}{2.8}{.75}{1}
      \Crossing{3}{2.8}{.75}{1}
    \fi
    \AnnularArc{0}{6.55}{7.75}
    \AnnularArc{.75}{5.85}{7}
    \AnnularArc{3}{5.15}{6.25}
    \AnnularArc{3.75}{4.45}{5.5}
    \node[axis] at (4.65,0) {$\times$};
  \end{tikzpicture}
}

\newcommand{\HalfDiagram}{%
  \begin{tikzpicture}[x=.7cm,y=.7cm,baseline=-.5ex]
    \FiveBlock{0}{1}
    \draw[strand] (-.8,2.15)..controls(-.15,2.15)and(0,2.15)..(0,2);
    \draw[strand] (-.8,-2.15)..controls(-.15,-2.15)and(0,-2.15)..(0,-2);
    \draw[strand] (.75,2)
      ..controls(.75,3.7)and(.9,4.05)..(2.2,4.05)
      ..controls(4.25,4.05)and(4.25,2.5)..(4.25,0)
      ..controls(4.25,-2.5)and(4.25,-4.05)..(2.2,-4.05)
      ..controls(.9,-4.05)and(.75,-3.7)..(.75,-2);
    \draw[strand] (1.5,2)
      ..controls(1.5,3.05)and(1.6,3.3)..(2.3,3.3)
      ..controls(3.5,3.3)and(3.5,2.5)..(3.5,0)
      ..controls(3.5,-2.5)and(3.5,-3.3)..(2.3,-3.3)
      ..controls(1.6,-3.3)and(1.5,-3.05)..(1.5,-2);
    \node[axis] at (2.55,0) {$\times$};
  \end{tikzpicture}%
}

\newcommand{\HalfState}[1]{
  \begin{tikzpicture}[x=.64cm,y=.64cm,baseline=-.5ex]
    \path[use as bounding box] (-.8,-2) rectangle (3.75,2);
    \ifnum#1=0
      \draw[strand] (-.75,.48)..controls(.05,.48)and(.05,-.48)..(-.75,-.48);
      \draw[strand] (2,0) ellipse (1.5 and 1.85);
      \draw[strand] (2,0) ellipse (.8 and 1.15);
    \else
      \draw[strand] (-.75,.48)
        ..controls(.25,.48)and(.05,1.85)..(1.6,1.85)
        ..controls(3.5,1.85)and(3.5,1.2)..(3.5,0)
        ..controls(3.5,-1.2)and(3.5,-1.85)..(1.6,-1.85)
        ..controls(.05,-1.85)and(.25,-.48)..(-.75,-.48);
      \draw[strand] (2,0) ellipse (.8 and 1.15);
    \fi
    \node[axis] at (2,0) {$\times$};
  \end{tikzpicture}%
}

\newcommand{\GuideResolution}[5]{
  \pgfmathsetmacro{\resolutionright}{#1+#3}
  \ifnum#4=1
    \draw[strand] (#1,#2+.34)
      ..controls(#1,#2+.04)and(\resolutionright,#2+.04)..(\resolutionright,#2+.34);
    \draw[strand] (#1,#2-.34)
      ..controls(#1,#2-.04)and(\resolutionright,#2-.04)..(\resolutionright,#2-.34);
    \draw[#5,line width=.65pt] (#1+.325,#2-.11)--(#1+.325,#2+.11);
  \else
    \draw[strand] (#1,#2-.34)--(#1,#2+.34);
    \draw[strand] (\resolutionright,#2-.34)--(\resolutionright,#2+.34);
    \draw[#5,line width=.65pt] (#1,#2)--(\resolutionright,#2);
  \fi
}

\newcommand{\GuideCut}[1]{
  \begin{tikzpicture}[x=.91cm,y=.78cm]
  \foreach \row/\leftnumber/\leftletter/\rightnumber/\rightletter in {
    0/1/-2/7/4,
    1/2/1/8/-5,
    2/3/1/9/-5,
    3/4/1/10/-5,
    4/5/-2/11/4,
    5/6/1/12/5}{
    \pgfmathsetmacro{\yc}{.36+.72*\row}
    \pgfmathsetmacro{\leftx}{.65*(abs(\leftletter)-1)}
    \pgfmathsetmacro{\leftxright}{\leftx+.65}
    \pgfmathsetmacro{\rightx}{.65*(abs(\rightletter)-1)}
    \pgfmathsetmacro{\rightxright}{\rightx+.65}
    \pgfmathtruncatemacro{\leftsign}{\leftletter>0 ? 1 : -1}
    \pgfmathtruncatemacro{\rightsign}{\rightletter>0 ? 1 : -1}
    \foreach \c in {0,1,2,3,4,5}{
      \pgfmathsetmacro{\xx}{.65*\c}
      \pgfmathtruncatemacro{\pass}{
        abs(\xx-\leftx)>.01 && abs(\xx-\leftxright)>.01 &&
        abs(\xx-\rightx)>.01 && abs(\xx-\rightxright)>.01}
      \ifnum\pass=1
        \draw[strand] (\xx,\yc-.36)--(\xx,\yc+.36);
      \else
        \draw[strand] (\xx,\yc-.36)--(\xx,\yc-.34);
        \draw[strand] (\xx,\yc+.34)--(\xx,\yc+.36);
      \fi
    }
    \ifnum#1=0
      \Crossing{\leftx}{\yc}{.65}{\leftsign}
      \Crossing{\rightx}{\yc}{.65}{\rightsign}
    \else
      \pgfmathtruncatemacro{\leftturnback}{\leftnumber==1 || \leftnumber==5}
      \def\rightcolour{red!70!black}
      \ifnum\rightnumber=8 \def\rightcolour{blue!70!black}\fi
      \ifnum\rightnumber=9 \def\rightcolour{blue!70!black}\fi
      \ifnum\rightnumber=10 \def\rightcolour{blue!70!black}\fi
      \GuideResolution{\leftx}{\yc}{.65}{\leftturnback}{red!70!black}
      \GuideResolution{\rightx}{\yc}{.65}{0}{\rightcolour}
    \fi
    \node[anchor=east,font=\scriptsize,text=black!70] at (-.18,\yc) {$\leftnumber$};
    \node[anchor=west,font=\scriptsize,text=black!70] at (3.43,\yc) {$\rightnumber$};
  }
  \draw[strand] (1.3,4.32)..controls(1.3,4.76)and(1.95,4.76)..(1.95,4.32);
  \draw[strand] (1.3,0)..controls(1.3,-.55)and(1.95,-.55)..(1.95,0);
  \foreach \c in {0,1,4,5}{
    \pgfmathsetmacro{\xx}{.65*\c}
    \draw[strand,->] (\xx,4.32)--(\xx,4.54);
    \draw[strand,->] (\xx,-.22)--(\xx,0);
    \ifnum#1=1 \node[font=\scriptsize,text=violet!80!black] at (\xx,-.56) {$v_+$};\fi
  }
  \ifnum#1=1
    \node[font=\scriptsize,text=violet!80!black,fill=white,inner sep=1pt] at (.975,1.80) {$w_-$};
  \fi
  \end{tikzpicture}
}

\title{A counterexample to the wrapping number conjecture}
\author{Qiuyu Ren}
\email{qren18@stanford.edu}
\hypersetup{
  pdftitle={A counterexample to the wrapping number conjecture},
  pdfauthor={Qiuyu Ren}
}

\begin{document}

\begin{abstract}
We exhibit an annular knot with wrapping number four whose Kauffman bracket has annular degree at most two. This disproves the wrapping number conjecture.
\end{abstract}

\maketitle

\section{Introduction}
For an annular link $L\subset S^1\times B^2$, its \emph{wrapping number} $\wrap(L)$ is the minimum of $|D\cap L|$ over meridional disks $D\subset S^1\times B^2$ transverse to $L$.

The Kauffman bracket skein module of $S^1\times B^2$, defined by the local relations
\[
\left\langle\SkeinCrossing\right\rangle
=A\left\langle\SkeinOne\right\rangle
+A^{-1}\left\langle\SkeinZero\right\rangle,
\qquad
\left\langle L\sqcup\SkeinCircle\right\rangle
=(-A^2-A^{-2})\langle L\rangle,
\]
is canonically isomorphic to $\Z[A^{\pm1}][z]$, where $z$ is represented by a core circle in $S^1\times B^2$ with the standard framing. Thus the bracket of a framed annular link is a polynomial $\langle L\rangle\in\Z[A^{\pm1}][z]$. Changing the framing multiplies it by a power of $-A^3$, so its \emph{annular degree} $\deg_z\langle L\rangle$ is an invariant of the underlying unoriented link. Hoste and Przytycki proposed that this degree always detects the wrapping number \cite{HP95}.

\begin{Conj}[Wrapping number conjecture {\cite[Conjecture~4]{HP95}} {\cite[Problem~1.28(a)]{K3}}]\label{conj:wrapping}
For every annular link $L$,
\[
\wrap(L)=\deg_z\langle L\rangle.
\]
\end{Conj}

Every resolution of a diagram meeting a radial seam $r$ times has at most $r$ essential circles. Hence $\deg_z\langle L\rangle\leq\wrap(L)$, and the content of the conjecture is the reverse inequality.

\begin{figure}[H]
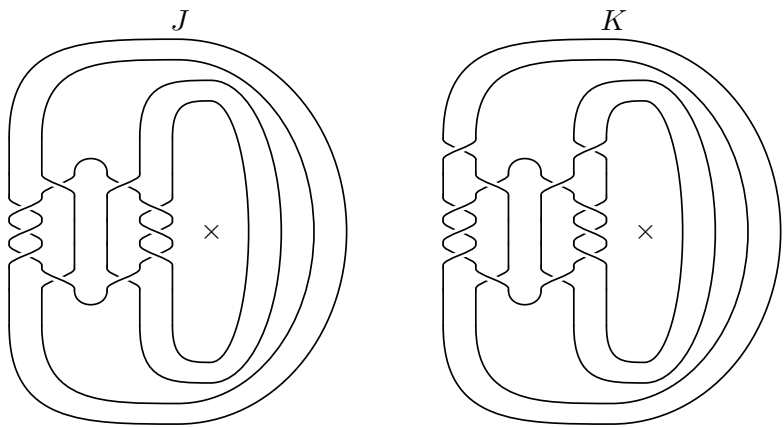

  \centering
  \begin{minipage}[t]{.28\linewidth}
    \centering
    $J$\par\smallskip
    \resizebox{\linewidth}{!}{\AnnularLink{0}}
  \end{minipage}\hspace{.06\linewidth}
  \begin{minipage}[t]{.28\linewidth}
    \centering
    $K$\par\smallskip
    \resizebox{\linewidth}{!}{\AnnularLink{1}}
  \end{minipage}
  \caption{The auxiliary annular link $J$ and the knot $K$. The mark $\times$ denotes the annular axis.}
  \label{fig:counterexample}
\end{figure}

\begin{Thm}\label{thm:counterexample}
The annular knot $K$ in Figure~\ref{fig:counterexample} satisfies
\[
\wrap(K)=4
\qquad\text{and}\qquad
\deg_z\langle K\rangle\leq2.
\]
In particular, Conjecture~\ref{conj:wrapping} is false.
\end{Thm}

Annular Khovanov homology categorifies the Kauffman bracket skein module of $S^1\times B^2$ \cite{APS04}. Let $\deg_k\operatorname{AKh}(L)=\max\{|k|:\operatorname{AKh}(L)_k\neq0\}$. In general, $\deg_z\langle L\rangle\leq\deg_k\operatorname{AKh}(L)\leq\wrap(L)$.

A weaker version of Conjecture~\ref{conj:wrapping}, known as the categorified wrapping number conjecture and proposed by Grigsby in 2010 \cite[Conjecture~1.2 and the discussion following it]{Martin23}, asserts that $\deg_k\operatorname{AKh}(L)=\wrap(L)$ for every annular link $L$. For the knot $K$, Section~\ref{sbsec:AKh} produces a nonzero class in $k$-grading $4$. Together with the upper bound visible in Figure~\ref{fig:counterexample}, this gives $\deg_k\operatorname{AKh}(K)=\wrap(K)=4$. Thus $K$ is not a counterexample to the categorified wrapping number conjecture.

In Section~\ref{sec:Kauffman}, we prove the assertion on the annular degree in two ways; the second argument also explains the motivation for the construction. In Section~\ref{sbsec:AKh}, we prove the assertion on the wrapping number using the annular Khovanov homology. In Section~\ref{sbsec:hyp}, we give an alternative proof of the assertion on the wrapping number using hyperbolic geometry. Similar arguments show that the $3$-component link $J$ in Figure~\ref{fig:counterexample} is also a counterexample to Conjecture~\ref{conj:wrapping}.

\begin{Rmk}
J\'ozef Przytycki suggested to us that $\wrap(K)=4$ may also be detected using the HOMFLY skein module, the two-variable Kauffman skein module (these skein modules of the solid torus were computed in \cite{Tur88} and \cite{Hoste1990}), or colored invariants. Indeed, a computer calculation shows that the Kauffman bracket of a $2$-cable of $K$ has annular degree $8$. Since this degree is at most $2\wrap(K)$, this recovers $\wrap(K)=4$.

It is therefore natural to ask, for example, whether the annular degrees of the colored Kauffman brackets detect the wrapping number. GPT-6 Astra suggested the following connection: even the assertion that an annular knot has wrapping number zero whenever all these degrees vanish would imply that the colored Jones polynomials detect the unknot. Indeed, suppose that a zero-framed knot $C\subset S^3$ has the same colored Jones polynomials as the unknot in every color, and let $P_C$ be one component of the untwisted Bing double of $C$, regarded as an annular knot in the complement of the other component. Since, over $\mathbb Q(A)$, the various Jones--Wenzl colorings of the core circle form a basis of the Kauffman bracket skein module of the solid torus, the colored Jones polynomials of $C$ determine those of its Bing double in all pairs of colors. Note also that the family of evaluations against all Jones--Wenzl colorings of the complementary core separates elements of the Kauffman bracket skein module of the solid torus. It follows that $P_C$ has the same colored annular Kauffman brackets as a local unknot, whereas $P_C$ has positive wrapping whenever $C$ is nontrivial.
\end{Rmk}

\subsection*{Statement on AI use}

The example $K$ was found autonomously by GPT-6 Astra, who essentially provided the proofs we present in Sections~\ref{sec:Kauffman} and~\ref{sbsec:AKh}. AI was used to draft this manuscript.

\subsection*{Acknowledgment}
We thank John Baldwin and J\'ozef Przytycki for helpful correspondence. This research was conducted during the period when the author served as a Clay Research Fellow.

\section{Annular degree}\label{sec:Kauffman}
The link $J$ in Figure~\ref{fig:counterexample} is obtained from $K$ by taking the braid-like resolution at each of its two uppermost crossings. At each of the two crossings added to obtain $K$, the braid-like resolution has coefficient $A$. Every other choice has lower annular degree. It follows that $\langle K\rangle=A^2\langle J\rangle+\cdots$, where $\cdots$ denotes terms with annular degree at most $2$. Thus, it suffices to prove that $\langle J\rangle$ has annular degree at most $2$.

\subsection{Method 1: Direct computation}
We cut $J$ into two mirror-image halves. The bracket of the right half is
\[
\begin{aligned}
\vcenter{\hbox{\scalebox{.42}{\HalfDiagram}}}
={}&A^{-1}\vcenter{\hbox{\scalebox{.48}{\HalfState{0}}}}
 +(A^5-A+2A^{-3}-A^{-7})\,
   \vcenter{\hbox{\scalebox{.48}{\HalfState{1}}}}\\
 &+\text{terms with zero wrapping}.
\end{aligned}
\]
The left half has the same computation except with $A$ replaced by $A^{-1}$. The terms with zero wrapping cannot contribute four essential circles after the two halves are glued together. Gluing the two displayed terms to their reflections therefore gives
\[
\begin{aligned}
\langle J\rangle
={}&\bigl((-A^2-A^{-2})
  +A(A^5-A+2A^{-3}-A^{-7})
  +A^{-1}(A^{-5}-A^{-1}+2A^3-A^7)\bigr)z^4+\cdots\\
={}&(-A^2-A^{-2}+A^2+A^{-2})z^4+\cdots
=\cdots.
\end{aligned}
\]
Here $\cdots$ denotes terms with annular degree at most $2$.

\subsection{Method 2: \texorpdfstring{$6$}{6}-strand Burau representation}
We sketch a less straightforward but more informative proof. We became aware of this argument only after further inquiry into the AI model, which revealed that it motivated the construction.

Let $V$ be the $\mathbb Q(A)$-vector space spanned by crossingless tangles with four lower and six upper endpoints, modulo the subspace spanned by those with fewer than four through strands. The latter cannot produce four essential circles after taking the annular closure. A basis of $V$ is $v_1,\ldots,v_5$, where $v_i$ has a cup joining the $i$th and $(i+1)$st upper endpoints and four through strands. The $6$-strand Temperley--Lieb algebra over $\mathbb Q(A)$, in which a circle evaluates to $-A^2-A^{-2}$, acts on $V$ by stacking. Let $E_i\in\operatorname{End}(V)$ be the action of the standard generator with a cap and a cup at positions $i$ and $i+1$ and vertical strands elsewhere. Then $E_iv_i=(-A^2-A^{-2})v_i$, $E_iv_{i\pm1}=v_i$, and $E_iv_j=0$ otherwise.

The bracket action of a positive braid generator $\sigma_i$ on $V$ is $AI+A^{-1}E_i$. Removing the scalar factor $A$ defines a representation $\rho\colon B_6\to\operatorname{GL}(V)$ by $\rho(\sigma_i)=I+A^{-2}E_i$. It is standard that $\rho$ is equivalent over $\mathbb Q(A)$ to the reduced Burau representation.

Let $\ell_3\colon V\to\mathbb Q(A)$ be the linear functional obtained by attaching a cap at positions $3$ and $4$ and extracting the coefficient of the identity matching on the remaining four strands. Directly from the diagrams, $E_3=v_3\ell_3$ is a rank-one operator.

The reduced Burau representation on $6$ strands is not faithful. We will use the following element Bigelow constructed in its kernel \cite[pp.~403--404]{Big99}. Namely, set
\[
\psi_1=\sigma_4\sigma_5^{-1}\sigma_2^{-1}\sigma_1,
\qquad
\psi_2=\sigma_4^{-1}\sigma_5^2\sigma_2\sigma_1^{-2}.
\]
Then, the commutator of $\psi_1^{-1}\sigma_3\psi_1$ and $\psi_2^{-1}\sigma_3\psi_2$ is in the kernel of $\rho$. Set $P_i=\rho(\psi_i)$ and $F_i=P_i^{-1}E_3P_i$. Expanding the resulting commutativity of $I+A^{-2}F_1$ and $I+A^{-2}F_2$ gives $F_1F_2=F_2F_1$. Each $F_i$ has image spanned by $P_i^{-1}v_3$. Using the action above, we obtain, with respect to the basis $v_1,\ldots,v_5$,
\[
\left.P_1^{-1}v_3\right|_{A=1}=(1,1,1,1,1)^T,
\qquad
\left.P_2^{-1}v_3\right|_{A=1}=(0,1,1,1,0)^T.
\]
Since these specializations are linearly independent, $P_1^{-1}v_3$ and $P_2^{-1}v_3$ span distinct lines over $\mathbb Q(A)$. The image of the common operator $F_1F_2=F_2F_1$ lies in their intersection, so $F_1F_2=0$. This gives
\[
0=F_1F_2=(P_1^{-1}v_3)\bigl(\ell_3P_1P_2^{-1}v_3\bigr)(\ell_3P_2).
\]
Since $P_1^{-1}v_3\ne0$ and $\ell_3P_2\ne0$, it follows that $\ell_3P_1P_2^{-1}v_3=0$.

Finally, commuting generators supported on the two disjoint triples of strands gives
\[
\psi_1\psi_2^{-1}
=(\sigma_2^{-1}\sigma_1^3\sigma_2^{-1})
 (\sigma_4\sigma_5^{-3}\sigma_4).
\]
The two factors on the right are precisely the two five-crossing blocks used to construct $J$. By the definitions of $V$ and $\ell_3$, the coefficient of $z^4$ in $\langle J\rangle$ is, up to the normalization of $\rho$, the matrix coefficient $\ell_3\rho(\psi_1\psi_2^{-1})v_3$. Since $\psi_1\psi_2^{-1}$ has exponent sum zero, no normalization factor appears. Thus, this coefficient is
\[
\ell_3\rho(\psi_1\psi_2^{-1})v_3
=\ell_3P_1P_2^{-1}v_3
=0.
\]

\section{Wrapping number}\label{sec:wrapping}

The diagram in Figure~\ref{fig:counterexample} meets a radial seam four times, so $\wrap(K)\leq4$. We give two proofs of the reverse inequality.

\subsection{Method 1: Annular Khovanov homology}\label{sbsec:AKh}

We use the unreduced integral annular Khovanov complex for annular link diagrams \cite{APS04}, with the conventions for labels and the annular grading from \cite[Section~2]{Rob13}. Although \cite{Rob13} works over $\mathbb F_2$, these conventions are coefficient-independent, and the local formula used below holds over $\Z$. A label $v_+$ or $v_-$ on an essential circle in a resolution of a diagram has annular degree $1$ or $-1$, respectively, while the labels $w_+$ and $w_-$ on a contractible circle have annular degree $0$. The subscripts $+$ and $-$ correspond to the standard Khovanov basis elements $1$ and $X$ of $\Z[X]/X^2$, respectively. The annular degree $k$ of an enhanced resolution (i.e.\ a fully labeled resolution) is the sum of these degrees, and the annular differential $d_0$ is the part of the ordinary Khovanov differential preserving $k$. We direct cube edges from the $0$-resolution to the $1$-resolution, where $0$ is the $A$-smoothing and $1$ is the $A^{-1}$-smoothing as in the Kauffman bracket skein relation.

\begin{figure}[tbp]
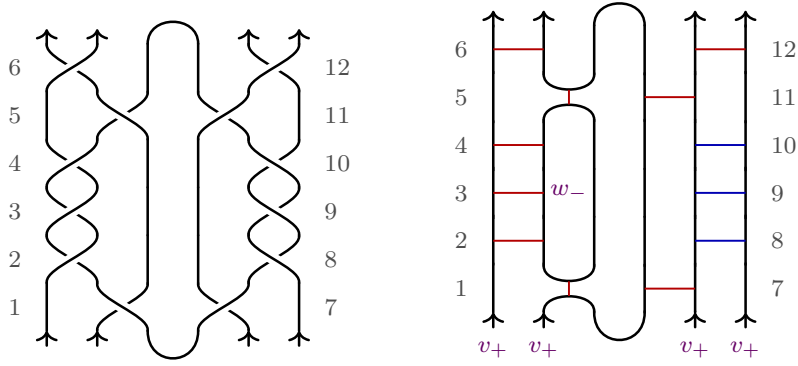

  \centering
  \begin{minipage}[t]{.30\linewidth}
    \centering
    \resizebox{\linewidth}{!}{\GuideCut{0}}
  \end{minipage}\hspace{.05\linewidth}
  \begin{minipage}[t]{.30\linewidth}
    \centering
    \resizebox{\linewidth}{!}{\GuideCut{1}}
  \end{minipage}
  \caption{The crossing numbering and the enhanced resolution $g$. The arrowheads indicate the annular identification, not an orientation. In the right panel, red and blue segments mark outgoing and incoming saddles, respectively.}
  \label{fig:akh-generator}
\end{figure}

Number the crossings of $K$ as in Figure~\ref{fig:akh-generator}, and consider the resolution corresponding to the string $s=000000011100$, where the leftmost digit corresponds to crossing $1$. The resolution consists of four essential circles and one contractible circle. Let $g$ be the corresponding element in the annular Khovanov chain complex obtained by labeling every essential circle with $v_+$ and the contractible circle with $w_-$. Then $k(g)=4$.

The only local formula we need is $m_0(v_+\otimes w_-)=0$. Indeed, the ordinary Frobenius multiplication sends $v_+\otimes w_-$ to $v_-$. This lowers the local annular degree, so its $k$-preserving part vanishes.

The outgoing cube edges from $s$ occur at crossings $1,2,3,4,5,6,7,11,12$. At crossings $6,7,11,12$, the target resolution has only two essential circles and hence no generator of annular degree $4$. At each of crossings $1,2,3,4,5$, the saddle merges the contractible circle, labeled $w_-$, with an essential circle labeled $v_+$; the resulting term vanishes by the local formula. Thus $d_0g=0$. Namely, $g$ is a cocycle in the annular Khovanov chain complex of $K$.

The only incoming edges come from changing crossing $8$, $9$, or $10$ from $1$ to $0$. Each source resolution has only two essential circles, so it has no generator of annular degree $4$. Consequently, $g$ is not a coboundary, so it represents a nonzero class in $\operatorname{AKh}(K)$ with $k$-degree $4$.

It is easy to see that $\deg_k\operatorname{AKh}(K)\leq\wrap(K)$: a diagram realizing the wrapping number meets a radial seam $r=\wrap(K)$ times, so every resolution has at most $r$ essential circles and every enhanced resolution has annular degree at most $r$ in absolute value. The nonzero class above therefore forces $\wrap(K)\geq4$.

\subsection{Method 2: Hyperbolic geometry}\label{sbsec:hyp}

We regard $S^1\times B^2$ as the standard solid torus in $S^3$, and let $U$ be the core of the complementary solid torus. Thus $U$ is the annular axis. The knot $K$ has winding number zero, so its wrapping number is even. If $\wrap(K)<4$, a meridional disk realizing the wrapping number meets $K$ in either zero or two points. After zero-surgery on $U$, this disk caps off to give an essential sphere or annulus in $S^1\times S^2\backslash K$. It therefore suffices to verify that this complement is hyperbolic.

This is certified by the following SnapPy \cite{SnapPy} computation in SageMath \cite{SageMath}.
\begin{Verbatim}[fontsize=\footnotesize,breaklines=true,breakanywhere=true]
sage: import snappy, sage.version
sage: sage.version.version, snappy.__version__
('10.7', '3.3.3a3')
sage: pd = [(25,5,26,4), (36,11,37,12), (10,35,11,36), (15,34,16,35), (28,34,29,33), (17,13,18,12), (14,22,15,21), (1,33,2,40), (31,27,32,26), (23,3,24,2), (3,25,4,24), (18,7,19,8), (20,9,21,10), (39,1,40,32), (5,30,6,31), (29,22,30,23), (8,19,9,20), (6,14,7,13), (37,16,38,17), (38,28,39,27)]
sage: M = snappy.Link(pd).exterior()
sage: M.dehn_fill((0,1), 1)
sage: N = M.filled_triangulation([1])
sage: N.verify_hyperbolicity()[0]
True
\end{Verbatim}
The PD code, obtained using KnotFolio \cite{KnotFolio}, lists $K$ first and $U$ second. Here $(0,1)$ is the preferred longitude of $U$. The last command uses interval arithmetic and therefore gives a rigorous verification \cite{HIKMOT}. Thus $\wrap(K)\geq4$.

\printbibliography

\end{document}